\documentclass[openany, a4paper, 11pt]{article}
\usepackage{graphicx}
\usepackage{mathrsfs}
\usepackage{amsfonts}

\usepackage{mathrsfs,amsfonts,amsmath}
\usepackage{color}
\usepackage{pstricks}
\usepackage{pst-plot}
\usepackage{pst-eps}
\usepackage{pst-grad}
 \newtheorem{theorem}{Theorem}[section]
 \newtheorem{definition}[theorem]{Definition}
 \newtheorem{lemma}[theorem]{Lemma}
 \newtheorem{corollary}[theorem]{Corollary}
 \newtheorem{proposition}[theorem]{Proposition}
 \newtheorem{remark}[theorem]{Remark}
 \newtheorem{condition}[theorem]{Condition}
 \newtheorem{example}{Example}[section]

 \def\blemma{\begin{lemma}}\def\elemma{\end{lemma}}
 \def\bproposition{\begin{proposition}}\def\eproposition{\end{proposition}}
 \def\btheorem{\begin{theorem}}\def\etheorem{\end{theorem}}
 \def\bcorollary{\begin{corollary}}\def\ecorollary{\end{corollary}}
 \def\bremark{\begin{remark}}\def\eremark{\end{remark}}
 \def\bcondition{\begin{condition}}\def\econdition{\end{condition}}

 \def\benumerate{\begin{enumerate}}\def\eenumerate{\end{enumerate}}
 \def\bitemize{\begin{itemize}}\def\eitemize{\end{itemize}}

 \def\beqlb{\begin{eqnarray}}\def\eeqlb{\end{eqnarray}}
 \def\beqnn{\begin{eqnarray*}}\def\eeqnn{\end{eqnarray*}}

 \def\ar{\!\!\!&}

 \def\mbb{\mathbb}

 \def\proof{\noindent{\it Proof.~~}}\def\qed{\hfill$\Box$\medskip}

 \def\bfP{{\mbox{\bf P}}}
 \def\bfE{{\mbox{\bf E}}}
 
 \def\bfX{{\mbox{\bf X}}}
 \def\bfY{{\mbox{\bf Y}}}
 \def\bfx{{\mbox{\bf x}}}
 \def\bfe{{\mbox{\bf e}}}
 \def\bfy{{\mbox{\bf y}}}
 \def\bfr{{\mbox{\bf r}}}

\begin{document}

\title{The peripatric coalescent}
\author{\textsc{By Amaury Lambert$^{1,2}$ and Chunhua Ma$^{1,2,3}$}
}
\date{}
\maketitle

\noindent\textsc{$^1$
UPMC Univ Paris 06\\
Laboratoire de Probabilit\'es et Mod\`eles Al\'eatoires CNRS UMR 7599}\\
\noindent\textsc{$^2$
Coll\`ege de France\\
Center for Interdisciplinary Research in Biology CNRS UMR 7241\\
Paris, France}\\\textsc{E-mail: }amaury.lambert@upmc.fr\\
\textsc{URL: }http://www.proba.jussieu.fr/pageperso/amaury/index.htm\\
\\
\noindent\textsc{$^3$Nankai University\\
School of Mathematical Sciences and LPMC \\
Tianjin, P.\,R.\ China}\\
\textsc{E-mail: }mach@nankai.edu.cn\\
\textsc{URL: } http://math.nankai.edu.cn/$\sim$mach

\begin{abstract}
\noindent
We consider a dynamic metapopulation involving one large population of size $N$ surrounded by colonies of size $\varepsilon_{_N}N$, usually called peripheral isolates in ecology, where $N\to\infty$ and $\varepsilon_{_N}\to 0$ in such a way that $\varepsilon_{_N}N\to\infty$.
The main population periodically sends propagules to found new colonies (emigration), and each colony eventually merges with the main population (fusion).
Our aim is to study the genealogical history of a finite number of lineages sampled at stationarity in such a metapopulation.

We make assumptions on model parameters ensuring that the total outer population has size of the order of $N$ and that each colony has a lifetime of the same order. We prove that under these assumptions, the scaling limit of the genealogical process of a finite sample is a censored coalescent where each lineage can be in one of two states: an inner lineage (belonging to the main population) or an outer lineage (belonging to some peripheral isolate). Lineages change state at constant rate and inner lineages (only) coalesce at constant rate per pair.

This two-state censored coalescent is also shown to converge weakly, as the landscape dynamics accelerate, to a time-changed Kingman coalescent.
\end{abstract}  	
\medskip
\textit{Running head.} The peripatric coalescent.

\section{Introduction}

\setcounter{equation}{0}

Many plant and animal populations in nature are highly fragmented, and this fragmentation plays a prominent role in the context of adaptation and speciation. Indeed, the emergence of new species is usually thought to be driven by geographical processes \cite{CO04}. First, \emph{allopatric speciation} occurs when various subpopulations belonging to the same initial species are separated by a geographical barrier that prevents hybridization between them (gene flow) and allows them to diverge (genetical differentiation) by local adaptation. Second, \emph{parapatric speciation} is a version of allopatric speciation where local adaptation is mediated by the existence of an environmental gradient (resource availability, environmental conditions). Third, when a species is present in one large, panmictic population surrounded by small colonies, usually called peripheral isolates, it is believed that the combination of founder events and of local adaptation to borderline environmental conditions leads to the formation of new species within the isolates. This phenomenon is called \emph{peripatric speciation}. We aim to study the genealogy of populations embedded in such a spatial context. The present study should then serve as a building brick for future work in the field of speciation modeling.

Population dynamic models specifying explicitly the spatial context are called \emph{metapopulation models} (Hanski and Gilpin \cite{HG97}). Typical such models  include: island model, isolation by distance, stepping stone models, extinction-recolonization models. From the point of view of speciation, all these models suffer from the same defect: they assume a given, constant number of subpopulations in the metapopulation, with fixed migration rates between them. As one of the authors of the present paper suggested (Lambert \cite{L10}), an alternative method would consist in considering species as ``spread out on a randomly evolving number of locations, allowing for repeated fragmentations of colonies, colonizations of new locations, as well as secondary contacts between subpopulations''. This author and others have designed such dynamic landscape models \cite{K00, ALC11, ACL13}, but usually in detailed ecological contexts whose study is only possible through numerical simulations (to the exception of \cite{ACL09}).

Here we propose a mathematical study of a dynamic landscape of the peripatric type. More specifically, we consider a dynamic population subdivision which involves one large main population surrounded by a random number of small peripheral isolates, that we will call colonies for simplicity. The size of the main population is constant equal to $N$, the size of each colony is constant equal to $\varepsilon_{_N}$ and the genealogy in each population is given by the Moran model. The number of colonies at time $t$ is denoted by $\xi_N(t)$.  The landscape dynamics is as follows (see Figure 1):
\begin{itemize}
\item At constant rate $\theta_N$, each individual sends independently $\varepsilon_{_N}N$ offspring to found a new colony;
\item
Each colony independently merges again with the main population at rate $\gamma_N \xi_N^{\alpha -1}$,
where $\alpha \ge 1$; at such so-called fusion time, $\varepsilon_{_N} N$ individuals among the new $(1+\varepsilon_{_N}) N$ individuals of the main population are chosen uniformly and simultaneously killed to keep its size constant.
\end{itemize}
Note that $(\xi_N(t);t\ge 0)$ is a pure-death process with immigration. The parameter $\alpha$ is meant to model the competition for space, since the fusion rate per colony grows with the number of colonies. This density-dependence disappears if $\alpha$ is chosen equal to 1.

The main purpose of this paper is to investigate the genealogy of a finite sample of lineages in the above peripatric metapopulation model. We will show that the history of such a sample, viewed backward in time, can be approximated, as $N\to\infty$ under certain assumptions, by a two-state censored coalescent, where the state of a lineage can be inner (lying in the main population) or outer (lying in a colony). Lineages change state at a constant rate per lineage, but only inner lineages can coalesce, at a constant rate per pair of lineages, as in Kingman coalescent \cite{K82}.

A two-state censored coalescent can be viewed as a new type of structured coalescent. The structured coalescent (see Takahata \cite{T88}, Notohara \cite{N90} and Herbots \cite{H97}) describes the ancestral genealogical process of a sample of lineages in a subdivided population connected by migration. The coalescent on two subpopulations was considered by \cite{T88}; for a finite number of subpopulations by \cite{N90}, and placed in a rigorous framework by \cite{H97}. To date, there have been a number of works dealing with the structured coalescent arising in various special types of metapopulations; see Nordborg and Krone \cite{NK02}, Eldon \cite{E09} and the references therein. Our results show that new types of structured coalescents can arise in some specific dynamic metapopulations.

We now give the heuristics giving rise to the result. We assume that $N\to\infty$ and $\varepsilon_{_N}\to 0$ in such a way that $\varepsilon_{_N}N\to\infty$, so that the size of colonies is large but neglectable compared to the main population (Assumption A). It is known that in a Moran model, inner lineages coalesce at constant rate per pair when time is rescaled by $N$ (Kingman coalescent \cite{K82}). We make assumptions on the parameters ensuring that all events changing the configuration of ancestral lineages occur on this time scale. This can only be done to the exception of coalescences in colonies, which happen instantaneously in the new time scale, leading to outer lineages which always all lie in different colonies. Also, in order to have a total outer population size of the order of $N$, we need to have a number of colonies of the order of $\varepsilon_{_N}^{-1}$. This can be achieved by the following choice of parameters (Assumption B). The \emph{per capita} emigration rate $\theta_N$ is taken equal to
$$
\theta_N= \frac{\theta}{\varepsilon_{_N} N^2} ,
$$
and the fusion rate $\gamma_N$ is taken equal to
$$
\gamma_N =  \gamma  \frac{\varepsilon_{_N}^{\alpha -1}}{N} .
$$
Under theses assumptions, the number of colonies is asymptotically deterministic, equal to $\varepsilon_{_N}^{-1}\,(\theta/\gamma)^{1/\alpha}$.

Now the rate at which a single inner lineage changes state is the rate at which a single lineage is taken in a fusion event (backward in time), which happens at rate
$$
\frac{\varepsilon_{_N}}{1+\varepsilon_{_N}} \,\gamma_N\, \xi_N^{\alpha} \approx \frac{\varepsilon_{_N}}{1+\varepsilon_{_N}} \,\gamma_N \,\varepsilon_{_N}^{-\alpha}\,\frac{\theta}{\gamma},
$$
which is equivalent to $\theta/N$ as $N\to\infty$. As a consequence, in the new time scale, inner lineages become outer lineages at contant rate $\theta$.

Also note that the probability that two lineages are taken in the same fusion vanishes, so that no two lineages can lie within the same colony. As a consequence, outer lineages are not allowed to coalesce.

Now the lifetime of a colony is approximately exponential with parameter
$$
\gamma_N\, \xi_N^{\alpha-1} \approx \gamma_N\,\varepsilon_{_N}^{1-\alpha}\,\left(\frac{\theta}{\gamma}\right)^{1-1/\alpha},
$$
which is equivalent to $(\theta/N)\,(\gamma/\theta)^{1/\alpha}$. As a consequence, in the new time scale, outer lineages become inner lineages at contant rate $\theta\,(\gamma/\theta)^{1/\alpha}$. By making these heuristics rigorous we get the results stated in Theorem \ref{thm}. Namely, the genealogical history of a finite sample of lineages, seen as a process backward in time, converges weakly (except at time 0, where instantaneous coalescences within colonies makes the limiting process not right-continuous) to the following two-state censored coalescent. Inner lineages coalesce at constant rate 1 per pair, and lineages change type at constant rate per lineage: inner lineages become outer lineages at rate $\theta$ and outer lineages become inner lineages at rate $\theta\,(\gamma/\theta)^{1/\alpha}$.

The paper is organized as follows. In Section 2, we give a detailed description of our
dynamic metapopulation model in forward and backward time. The main result, Theorem \ref{thm}, is
stated in Section 3. In addition, we also prove that under fast landscape dynamics, the censored coalescent
converges weakly to a time-changed version of the Kingman coalescent \cite{K82}. Finally, a section is dedicated to the formal proofs of the above results.

\section{Metapopulation model}
\setcounter{equation}{0}

\subsection{Forward dynamics}

Let $N\in\mathbb{N}$ with $\mathbb{N}:=\{0,1,2,\cdots\}$ and
let $\varepsilon_{_N}$ be any positive number such that $\varepsilon_{_N} N\in\mathbb{N}$.
Let $\theta_N$, $\gamma_N$ and $\alpha$
be positive constants.
Consider a dynamic metapopulation model involving one large population of size $N$, called main population, and a random number of small populations, called colonies, of size $\varepsilon_{_N}N$. The main population periodically sends propagules (or emigrants) that found new colonies and ultimately each colony merges again with the main population. A further assumption is as follows.
See Figure 1 for an illustration.

\noindent{\rm (1)} The number of colonies, denoted by $\{\xi_N(t): t\geq0\}$, evolves as
a pure death density-dependent process with immigration and the transition rates are given by
\beqlb \begin{array}{lll}
     j\rightarrow j+1 &\text{ at rate } & N\theta_N,\\
     j\rightarrow j-1 &\text{ at rate } & \gamma_Nj^\alpha.
    \end{array}
   \label{transition rate}
  \eeqlb
When $\alpha=1$, the process $\{\xi_N(t)\}$ is reduced to a pure death branching process with immigration.
It follows from Kelly \cite{K79} that $\{\xi_N(t)\}$
with any initial value has the stationary distribution $\pi_N$ given by
 \beqlb
 \pi_N(0)=\Big(1+\sum_{j=1}^\infty\frac{(N\theta_N/\gamma_N)^j}{(j!)^\alpha}\Big)^{-1}
 \ \mbox{and}\
 \pi_N(k)=\frac{(N\theta_N/\gamma_N)^k}{(k!)^\alpha}
 \Big(1+\sum_{j=1}^\infty\frac{(N\theta_N/\gamma_N)^j}{(j!)^\alpha}\Big)^{-1}
  \label{stationary distribution}
 \eeqlb
for $k\geq1$. We assume that $\xi_N(0)$ is distributed as $\pi_N$. Then $\{\xi_N(t)\}$ is a
stationary Markov chain. Let $(P_t^N)_{t\geq0}$ be its semigroup. For any finite set $\{t_1<
t_2< \cdots< t_n\}\subset \mbb{R}$ define the probability measure
on $\mathbb{N}$ by
 \beqlb\label{2.16}
\eta^N_{t_1,t_2,\cdots,t_n}(j_1,j_2,\cdots,j_n)
 =
\pi^N(j_1)P^N_{t_2-t_1}(j_1,j_2)\cdots P^N_{t_n-t_{n-1}}(j_{n-1},j_n).
 \eeqlb
Then $\{\eta^N_{t_1,t_2,\cdots,t_n}: t_1< t_2< \cdots< t_n\in \mbb{R}\}$ is
a consistent family. By Kolmogorov's theorem, there is a stochastic
process $\{\xi_N(t): t\in \mbb{R}\}$ with finite-dimensional distributions
given by (\ref{2.16}). Clearly, $\{\xi_N(t): t\in \mbb{R}\}$ is a stationary
Markov chain with one-dimensional marginal distribution $\pi_N$ and
transition semigroup $(P^N_t)_{t\ge 0}$.

\noindent{\rm (2)} At the jump times of $\xi_N(t)$ from $j$ to $j+1$, one individual, chosen uniformly
at random from the large population, gives birth to $\varepsilon_{_N}N$ emigrant offspring individuals which
found a new colony. We refer to such an event as ``emigration'' (of new colonies) or ``fission".

\noindent{\rm (3)} At the jump times of $\xi_N(t)$ from $j$ to $j-1$, one colony is chosen at random
from the $j$ current colonies and all the $\varepsilon_{_N}N$ individuals within this colony immediately migrate back into the main population. We refer to such an event as a ``fusion'' (of colonies with the main population). Instead of keeping all those $(1+\varepsilon_{_N})N$ individuals in the main population alive, only $N$ of them survive this fusion event, which are chosen uniformly at random
among the $(1+\varepsilon_{_N})N$ previously existing individuals.

\noindent{\rm (4)} Between the jump times of $\xi_N(t)$, the large population and the colonies independently
evolve as Moran models, that is, at rate $1$ each individual independently gives birth to
a single offspring, and simultaneously a uniformly chosen individual is killed.

\bigskip

\includegraphics[width=0.95\textwidth]{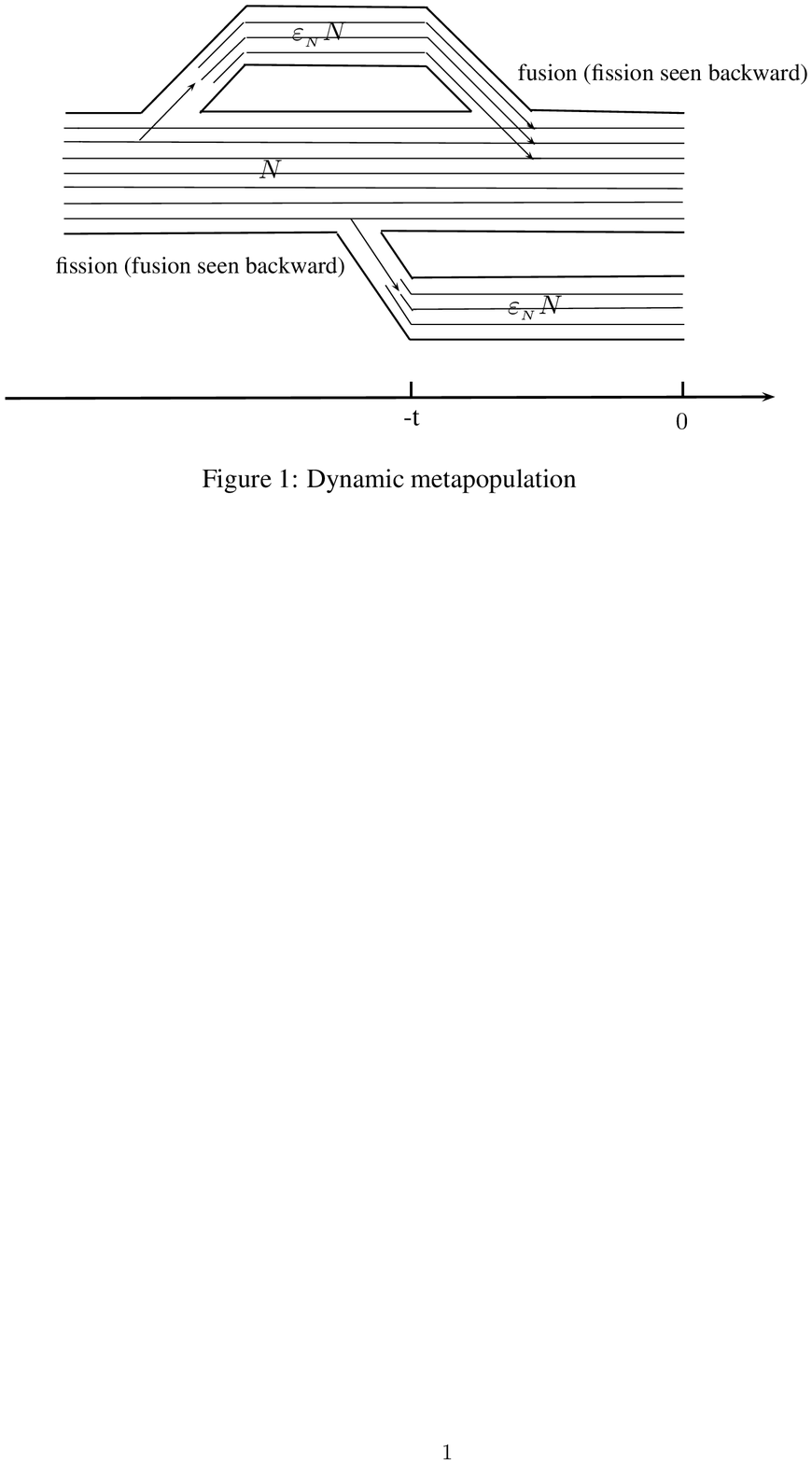}

\subsection{Backward dynamics}

Now we start with a sample of $n$ lineages at time 0 and proceed backward in time. Let $\bfX_N(t)=(X_N^0(t),X_N^1(t),\cdots,X_N^n(t))$  be the ancestral process of this sample defined for $t\ge 0$ by
\begin{itemize}
\item[\;]$X_N^0(t)=$ the number of lineages in the main population at time $-t$,

\item[\;]${X}_N^i(t)=$ the number of colonies containing $i$ lineages
at time $-t$ \ ($1\leq i\leq n$).
\end{itemize}
We set $\bfX_N(0)=\bfx$, where $\bfx=(x_0,x_1,\cdots,x_n)\in\mathbb{N}^{n+1}$ with $x_0+\sum_{j=1}^njx_j=n$.
The process $\{\bfX_N(t): t\geq0 \}$ has state-space
 \beqnn
E:=\Big\{(x_0,x_1,\cdots,x_{n})\in
\mathbb{N}^{n+1}: 1\leq x_0+\sum_{j=1}^{n}jx_j\leq n \Big\}.
 \eeqnn
Define the subspace $\Pi$ of $E$ by
\beqnn
 \Pi:=\big\{(x_0,x_1,0,\cdots,0)\in
\mathbb{N}^{n+1}: 1\leq x_0+x_1\leq n\big\}.
 \eeqnn
Consider the projection $\Gamma: (x_0,x_1,0,\cdots,0)\mapsto (x_0,x_1)$ from $\Pi$ to $\mathbb{N}^2$.
\beqnn
\Gamma(\Pi)=\big\{(x_0,x_1)\in
\mathbb{N}^{2}: 1\leq x_0+x_1\leq n\big\}.
 \eeqnn
By the action of the homeomorphism $\Gamma$, $\Gamma(\Pi)$ can be regarded as a subspace of $E$, and we thus still denote it by $\Pi$ for simplicity.
For $\bfx\in E$, let
$$
\bar{{\bf x}}:=\Big(x_0,\sum_{j=1}^nx_j\Big).
$$
We will use this notation for the following reason. Because in the new time scale lineages lying in the same colony immediately coalesce, the configuration ${\bf x}$ immediately turns into $(x_0,\sum_{j=1}^nx_j, 0,\ldots, 0)$ where all outer lineages are now alone in their respective colonies.
Note that $\bfx\mapsto\bar{\bfx}$ is an injection from $E$ to $\Pi$.
We also write
$\bfe_j=(0,\ldots,0,1,0,\ldots,0)\in\mathbb{N}^{n+1}$ whose $(j+1)$-th component is $1$ for
$j=0,\ldots,n$.

Let $\eta_N(t)=\xi_N(-t)$ for $t\geq0$. It follows from \cite[Lemma 1.5, P.9]{K79} that
$\{\eta_N(t): t\geq0 \}$ is still a stationary Markov process
with the same transition rates as (\ref{transition rate}). Thus,
the fission events (fusions seen backward) happen at rate $\theta_NN$ and, conditioned on $\eta_N(t)$, the fusion
events (fissions seen backward) happen at rate $\gamma_N\eta_N^{\alpha}(t)$.
At any fission time, every lineage independently exits from the main population with probability $\varepsilon_{_N}
/(1+\varepsilon_{_N})$. At any fusion time, one colony is chosen at random from the existing colonies
and the (say) $i$ lineages in this colony enter the main population, and simultaneously
coalesce together (if $i\ge 2$), and coalesce with their ancestor in the main population (if it is also in the sample; but asymptotically, with high probability $i=1$ and the ancestor is not in the sample). Between
fission and fusion times, coalescences within the main population or within colonies may happen.
We again refer to Figure 1 for an illustration.

Based on the above description, it is not hard
to see that $\{(\bfX_N(t),\eta_N(t)): t\geq0\}$ is a time-homogeneous Markov chain taking values in
 $E\times \mathbb{R}_+$. The corresponding generator is given by
  \beqlb
 \bar{A}_Ng(\bfx,k)=\bar{\psi}_Ng(\bfx,k)+\bar{\phi}_Ng(\bfx,k)+\bar{\Gamma}_Ng(\bfx,k)
 \label{bar A_N}
  \eeqlb
 for any bounded function $g$ on $E\times\mathbb{N}$. Here
  \beqnn
  \bar{\psi}_Ng(\bfx,k)=\sum_{j=2}^{n}x_j\binom{j}{2}\frac{2}{\varepsilon_{_N}N-1}
  \Big(g(\bfx-\bfe_j+\bfe_{j-1},k)-g(\bfx,k)\Big),
  \eeqnn
 which corresponds to coalescences within colonies. Note that $\psi_Ng(\bfx,u)\equiv0$ if
 $\bfx\in\Pi$. Then
  \beqnn
  \bar{\phi}_Ng(\bfx,k)
  \ar=\ar
  \binom{x_0}{2}\frac{2}{N-1}\Big(g(\bfx-\bfe_1,k)-g(\bfx,k)\Big)\\
  \ar \ar
  +\,N\theta_{N} \sum_{r=1}^{x_0}\binom{x_0}{r}\Big(\frac{\varepsilon_{_N}}{1+\varepsilon_{_N}}\Big)^r
  \Big(\frac{1}{1+\varepsilon_{_N}}\Big)^{x_0-r}
  \Big(g(\bfx-r\bfe_0+\bfe_{r},k+1)-g(\bfx,k)\Big)\\
  \ar \ar
  +\,\gamma_N k^{\alpha}(x_1/k)(1-(x_0/N))
  \Big(g(\bfx-\bfe_1+\bfe_0,k-1)-g(\bfx,k)\Big)1_{\{k>0\}}\\
  \ar \ar
  +\,\gamma_Nk^{\alpha}\sum_{j=2}^{N}(x_j/k)(1-(x_0/N))
  \Big(g(\bfx-\bfe_j+\bfe_0,k-1)-g(\bfx,k)\Big)1_{\{k>0\}}\\
  \ar \ar
  +\,\gamma_Nk^{\alpha}\sum_{j=1}^{N}(x_j/k)(x_0/N)
  \Big(g(\bfx-\bfe_j,k-1)-g(\bfx,k)\Big)1_{\{k>0\}}.
  \eeqnn
In $\bar{\phi}_N$, the first term corresponds to coalescences within the main population, the second term corresponds to exit from the main population, the third term corresponds to entrance into the main population, the last two terms correspond to simultaneous entrance into the main population and coalescence. The fourth term is identically equal to $0$ if $\bfx\in\Pi$. Last,
  \beqnn
  \bar{\Gamma}_Ng(\bfx,k)
  \ar=\ar
  N\theta_N\Big(\frac{1}{1+\varepsilon_{_N}}\Big)^{x_0}
  \Big(g(\bfx,k+1)-g(\bfx,k)\Big)\\
  \ar \ar\,
  +\gamma_Nk^{\alpha}\Big(1-\sum_{j=1}^{N}x_j/k\Big)
  \Big(g(\bfx,k-1)-g(\bfx,k)\Big)1_{\{k>0\}},
  \eeqnn
 which corresponds to the event that the number of colonies increases or decreases
 but the ancestral process does not change.

\section{Convergence to the two-state censored coalescent}

\setcounter{equation}{0}

\subsection{Main results}

Let $D([0,\infty), S)$
be the space of all c\`{a}dl\`{a}g functions $x: [0,\infty) \rightarrow S$ endowed with the
Skorokhod topology for any separable and complete metric space $S$;
see Ethier and Kurtz \cite[p.116]{EK86} for details.
For $N\in\mathbb{N}$, we consider the sequence of processes
$\{(\bfX_N(\cdot), \eta_N(\cdot))\}$. Define
 \beqnn
\bfY_N(t)=\bfX_N(Nt)\quad\mbox{and}\quad\tilde{\eta}_N(t)=\varepsilon_{_N}\eta_N(Nt).
 \eeqnn
Let $\theta>0$ and $\gamma>0$ be constants.
We further assume the following conditions:
\begin{itemize}

 \item[(A)]\;$\varepsilon=\varepsilon_{_N}$ satisfying $\varepsilon_{_N}\rightarrow0$
 and $\varepsilon_{_N} N\rightarrow\infty$ as $N\rightarrow\infty$;

 \item[(B)]\;$\theta_N=\theta/(\varepsilon_N N^2)$ and $\gamma_N=\gamma\varepsilon_{_N}^{\alpha-1}/N$.
\end{itemize}
Recall that $ \bfy\in E$ and the corresponding $\bar{\bfy}\in\Pi$. The main result of the paper follows.
\btheorem\label{thm}\;Let $\bfx\in E$. Under conditions (A) and (B),
the finite-dimensional distributions of the ancestral process
$\{\bfY_N(t),\,t\geq0\}$ starting at $\bfy$ converges to
those of a $\Pi$-valued continuous time Markov chain $\{\bfY(t),\ t\geq0\}$ starting at
$\bar{\bfy}$, except at time $0$.
The corresponding infinitesimal generator ${\bf Q}=(q_{\bfr,\bfr'})_{\bfr,\bfr'\in\Pi}$ is given by
  \beqlb q_{\bfr,\bfr'} = \left\{ \begin{array}{ll}
     -(\theta r_0+\theta(\gamma/\theta)^{1/\alpha}r_1+r_0(r_0-1)), &
     \text{ if } \bfr'=\bfr,\\
     \theta r_0, &
     \text{ if } r_0\neq 0\text{ and }\bfr'=\bfr+(-1,1),  \\
     \theta(\gamma/\theta)^{1/\alpha}r_1, &
     \text{ if } r_1\neq0\text{ and } \bfr'=\bfr+(1,-1),\\
     r_0(r_0-1), &
     \text{ if } \bfr'=\bfr+(-1,0),\\
     0, &
     \text{ otherwise.}
   \end{array} \right.
   \label{Y}
  \eeqlb
where $\bfr=(r_0,r_1)\in\Pi$. Furthermore, if the initial value $\bfy\in\Pi$,
weak convergence on $D([0,\infty), \Pi)$ to $\{\bfY(t)\}$ holds.
\etheorem
The previous statement describes the asymptotic genealogical history of a finite sample of lineages, seen as a process backward in time. Except at time 0, where instantaneous coalescences within colonies makes the limiting process not right-continuous, this process converges weakly to a two-state censored coalescent, where type 0 corresponds to inner lineages (lying in the main population) and type 1 to outer lineages (lying in pairwise distinct colonies). Inner lineages coalesce at constant rate 1 per (ordered) pair, and lineages change type at constant rate per lineage: inner lineages become outer lineages at rate $\theta$ and outer lineages become inner lineages at rate $\theta\,(\gamma/\theta)^{1/\alpha}$.\\
\\
Now consider a sequence of censored coalescent processes $\{\bfY_k(t)\}$ defined by (\ref{Y}) with parameters $\theta$ and $\gamma$
replaced by $\theta_k$ and $\gamma_k$, and the initial value $\bfY_k(0)=\bfy\in\Pi$ with
$y_0+y_1=n$.
Let $Y_k(t)=Y_k^0(t)+Y_k^1(t)$ and let $I_n=\{0,1,2,\cdots,n\}$.  We assume that
 \begin{itemize}
 \item[(C)]\;As $k\rightarrow\infty$, $\theta_k\rightarrow\infty$, $\gamma_k\rightarrow\infty$ and
 $\theta_k/\gamma_k\rightarrow p$ for some constant $p>0$.
 \end{itemize}
The above condition corresponds to the acceleration of the landscape dynamics (emigration and fusion). The following theorem states that such an acceleration gives rise to a single state coalescent process, where coalescence rates are obtained by averaging over the probability of presence in the main population.

\btheorem\label{thm2}\; Under condition (C), the process $\{Y_k(t),\,t\geq0\}$ starting at $n$ converges weakly
to the time-changed  $n$-Kingman coalescent process $\{K(t),\,t\geq0\}$ on $D([0,\infty), I_n)$.
When $K=l$, the coalescence rate is given by
 \beqnn
 c_l=\sum_{j=1}^lj(j-1)\binom{l}{j}\Big(\frac{p^{1/\alpha}}{1+p^{1/\alpha}}\Big)^j
 \Big(\frac{1}{1+p^{1/\alpha}}\Big)^{l-j}.
  \eeqnn
 \etheorem
\bremark It is easy to see that if $p=0$ which corresponds to predominant emigrations,
$\{Y_k(t),\,t\geq0\}$ converges weakly
to the constant process $\{K(t)\equiv n, t\geq0\}$; if $p=\infty$ which corresponds to predominant fusions,
$\{Y_k(t),\,t\geq0\}$ converges weakly
to the standard Kingman coalescent $\{K(t), t\geq0\}$ (i.e., $c_l=l(l-1)$).
\eremark
\subsection{Proofs}

To prove Theorem \ref{thm}, we start by proving the following lemmas.

\blemma\label{lemma 1}\;Under conditions (A) and (B), as $N\rightarrow\infty$,
  \beqnn
  \tilde{\eta}_N(\cdot)\overset{p}{\longrightarrow}(\theta/\gamma)^{1/\alpha}
  \eeqnn
in $D([0,\infty), \mathbb{R}_+)$.
\elemma

\proof Recall that the number  $\xi_N(Nt)$ of colonies of size $\varepsilon_{_N}N$ is a pure death density-dependent
 process with immigration with  transition rates given by
\beqnn \begin{array}{lll}
     j\rightarrow j+1 &\text{ at rate } & \theta/\varepsilon_{_N},\\
     j\rightarrow j-1 &\text{ at rate } & \gamma\varepsilon_{_N}^{\alpha-1}j^\alpha.
    \end{array}
     \eeqnn
By (\ref{stationary distribution}), it has the stationary distribution $\pi_N$ given by
 \beqnn
 \pi_N(0)=\Big(1+\sum_{j=1}^\infty\frac{(\theta/(\gamma\varepsilon_{_N}^\alpha))^j}{(j!)^\alpha}\Big)^{-1}
 \ \mbox{and}\
 \pi_N(k)=\frac{(\theta/(\gamma\varepsilon_{_N}^\alpha))^k}{( k!)^\alpha}
 \Big(1+\sum_{j=1}^\infty\frac{(\theta/(\gamma\varepsilon_{_N}^\alpha))^j}{(j!)^\alpha}\Big)^{-1}.
 \eeqnn
Recall that $\xi_N(0)$ is distributed as $\pi_N$.
Let $M_{N}=[(\theta/\gamma)^{\frac{1}{\alpha}}/\varepsilon_{_N}]$.
For $k>M_N$,
 \beqnn
 \ar\ar\pi_N(k)\\
 \ar=\ar\frac{(\theta/(\gamma\varepsilon_{_N}^\alpha))^{M_N}(\theta/(\gamma\varepsilon_{_N}^\alpha))^{k-M_N}}
{(M_N!\prod_{j=M_N+1}^k j)^\alpha}
  \bigg(1+\Big(\sum_{j=1}^{M_N}+\sum_{j=M_N+1}^\infty\Big)\frac{(\theta/(\gamma\varepsilon_{_N}^\alpha))^j}{(j!)^\alpha}\bigg)^{-1}\\
  \ar=\ar
  \prod_{j=M_N+1}^k\frac{\theta/(\gamma\varepsilon_{_N}^\alpha)}{j^\alpha}
  \bigg(\frac{(M_N!)^\alpha}{(\theta/(\gamma\varepsilon_{_N}^\alpha))^{M_N}}
  +\sum_{j=1}^{M_N-1}\prod_{i=j+1}^{M_N}\frac{i^\alpha}{\theta/(\gamma\varepsilon_{_N}^\alpha)}
  +1
  +\sum_{j=M_N+1}^\infty\prod_{i=M_N+1}^j\frac{\theta/(\gamma\varepsilon_{_N}^\alpha)}{i^\alpha}\bigg)^{-1}.
 \eeqnn
Note that if $j>M_N$, $j^\alpha>\theta/(\gamma\varepsilon_{_N}^\alpha)$ and if $j>3M_N$,
$\frac{\theta/(\gamma\varepsilon_{_N}^\alpha)}{j^\alpha}\leq
\frac{\theta/(\gamma\varepsilon_{_N}^\alpha)}{(3M_N)^\alpha}<\frac{1}{2^\alpha}$ for sufficiently large $N$.
Then $\prod_{j=M_N+1}^{4M_N}\frac{\theta/(\gamma\varepsilon_{_N}^\alpha)}{j^\alpha}\leq 2^{-\alpha M_N}$ and
  \beqnn
\sum_{k=4M_N+1}^\infty\prod_{j=M_N+1}^k\frac{\theta/(\gamma\varepsilon_{_N}^\alpha)}{j^\alpha}
 \leq
\prod_{j=M_N+1}^{4M_N}\frac{\theta/(\gamma\varepsilon_{_N}^\alpha)}{j^\alpha}\sum_{k=1}^\infty2^{-\alpha k}\leq O(2^{-\alpha M_N}).
  \eeqnn
If follows that $\pi_N([4M_N,\infty))\leq O(2^{-\alpha M_N})$. Thus the sequence $\{\varepsilon_{_N}\xi_N(0)\}$ is tight.
On the other hand, $\{\varepsilon_{_N}\xi_N(Nt)\}$  takes values in $\{i\varepsilon_{_N}: i\in\mathbb{N}\}$ and its generator is given by
 \beqnn
 L_Nf(z)=\gamma\varepsilon_{_N}^{\alpha-1}(z/\varepsilon_{_N})^\alpha(f(z-\varepsilon_{_N})-f(z))
 +(\theta/\varepsilon_{_N})(f(z+\varepsilon_{_N})-f(z)),
 \eeqnn
for any continuous bounded function $f$ on $\mathbb{R}_+$. Let $C_c^2(\mathbb{R}_+)$ be the set of
 twice differentiable functions with compact support on $\mathbb{R}_+$.
 It is not hard to see that as $N\rightarrow\infty$,
for $f\in C_c^2(\mathbb{R}_+)$,
 \beqlb
 \|L_Nf-Lf\|\rightarrow0 \mbox{ and } Lf(z)=(\theta-\gamma z^\alpha)f'(z),
 \label{xi}
 \eeqlb
where $\|f\|=\sup_{x\in\mathbb{R}_+}|f(x)|$. The Markov process $\xi$ with generator $L$ is actually deterministic and satisfies the ODE:
 \beqnn
 \xi'(t)=\theta-\gamma\xi^\alpha(t),
 \eeqnn
which has the unique equilibrium point $(\theta/\gamma)^{1/\alpha}$. It follows from (\ref{xi}), \cite[Theorem 6.1, P.28]{EK86}
and \cite[Theorem 9.10, P.244]{EK86} that $\varepsilon_{_N}\xi_N(0)\overset{w}{\longrightarrow} (\theta/\gamma)^{1/\alpha}$
as $N\rightarrow\infty$. Again by (\ref{xi}), \cite[Corollary 8.7, p.231]{EK86} shows that
$\{\xi_N(t):t\geq0\}$ converges weakly to the constant function $\{\xi(t)\equiv(\theta/\gamma)^{1/\alpha},\ t\geq0\}$
on $D([0,\infty), \mathbb{R}_+)$. Since $\xi_N(\cdot)$ is stationary and time-reversible, the same weak convergence holds
for $\tilde{\eta}_N(\cdot)$. Then this lemma is proved. \qed

 As in Section 2 it is easy to see that
 $(\bfY_N(\cdot), \tilde{\eta}_N(\cdot))$ is a continuous time Markov chain taking values in
 $E\times\mathbb{R}_+$. Based on (\ref{bar A_N}) and Conditions (A) and (B), a simple calculation shows that
 the corresponding generator is given by
  \beqlb
 A_Ng(\bfy,u)=\psi_Ng(\bfy,u)+\phi_Ng(\bfy,u)+\Gamma_Ng(\bfy,u)
 \label{A_N}
  \eeqlb
 for any bounded function $g$ on $\mathbb{R}_+\times E$. Here
  \beqnn
  \psi_Ng(\bfy,u)=2\sum_{j=2}^{n}y_j\binom{j}{2}\frac{1}{\varepsilon_{_N}}
  \Big(1+\frac{1}{\varepsilon_{_N}N-1}\Big)
  \Big(g(\bfy-\bfe_j+\bfe_{j-1},u)-g(\bfy,u)\Big).
  \eeqnn
 Note that $1/(\varepsilon_{_N}N-1)\rightarrow0$ as $N\rightarrow\infty$ by Condition (A). We also have
  \beqnn
  \phi_Ng(\bfy,u)
  \ar=\ar
  2\binom{y_0}{2}\Big(g(\bfy-\bfe_0,u)-g(\bfy,u)\Big)\\
  \ar \ar
  +\,\theta y_0
  \Big(g(\bfy-\bfe_0+\bfe_1,u+\varepsilon_{_N})-g(\bfy,u)\Big)\\
  \ar \ar
  +\,\gamma u^{\alpha-1}y_1
  \Big(g(\bfy-\bfe_1+\bfe_0,u-\varepsilon_{_N})-g(\bfy,u)\Big)1_{\{u>0\}}\\
  \ar \ar
  +\,\gamma u^{\alpha-1}\sum_{j=2}^{N}y_j
  \Big(g(\bfy-\bfe_j+\bfe_0,u-\varepsilon_{_N})-g(\bfy,u)\Big)1_{\{u>0\}}\\
  \ar \ar
  +\,\Big(\varepsilon_{_N}R_{1,N}g(\bfy,u)+\frac{1}{N}u^{\alpha-1}1_{\{u>0\}}
  R_{2,N}g(\bfy,u)\Big).
    \eeqnn
Here the fourth term is identically equal to $0$ if $\bfy\in\Pi$. In the last term, $R_{1,N}$ and $R_{2,N}$
are bounded linear operators satisfying $\|R_{i,N}\|\leq C$ for some
constant $C$. This last term includes the simultaneous entrance into the main population and coalescence
and it vanishes if $c_1\leq u\leq c_2$ for positive numbers $c_1$ and $c_2$. Last, we have
  \beqnn
  \Gamma_Ng(\bfy,u)
  \ar=\ar
  \theta\varepsilon_{_N}^{-1}(1-y_0\varepsilon_{_N})
  \Big(g(\bfy,u+\varepsilon_{_N})-g(\bfy,u)\Big)\\
  \ar \ar\,
  +\gamma u^{\alpha}\varepsilon_{_N}^{-1}\Big(1-\varepsilon_{_N}u^{-1}\sum_{j=1}^Ny_j\Big)
  \Big(g(\bfy,u-\varepsilon_{_N})-g(\bfy,u)\Big)1_{\{u>0\}}.
  \eeqnn
Let us write $c_{\psi}^N(\bfy)$ (resp. $c_{\phi}^N(\bfy,u), c_{\Gamma}^N(\bfy,u))$
the total rate of the events generated by $\psi_N$ (resp. $\phi_N$, $\Gamma_N$) when $A_Ng$
is applied to $(\bfy,u)$.
Then
 \beqnn
 c_{\psi}^N(\bfy)
  =2\sum_{j=2}^{n}y_j\binom{j}{2}\frac{1}{\varepsilon_{_N}}
  \Big(1+\frac{1}{\varepsilon_{_N}N-1}\Big),
 \eeqnn
 \beqnn
 c_{\phi}^N(\bfy,u)
  =
 2\binom{y_0}{2}+\theta y_0
 +\gamma u^{\alpha-1}1_{\{u>0\}}\sum_{j=1}^{N}y_j+\varepsilon_{_N}(1+u^{\alpha-1}1_{\{u>0\}})
 \eeqnn
 and
 \beqnn
 c_{\Gamma}^N(\bfy,u)
 =
 \theta\varepsilon_{_N}^{-1}(1-y_1\varepsilon_{_N})
 +\gamma u^{\alpha}\varepsilon_{_N}^{-1}\Big(1-\varepsilon_{_N}u^{-1}\sum_{j=1}^Ny_j\Big)1_{\{u>0\}}.
 \eeqnn
 Let us introduce the following notation,
  \beqnn
 \sigma^N_0=\inf\{t\geq0: \bfY_N(t)\in\Pi\}
  \eeqnn
 and
  \beqnn
  \sigma_1^N=\inf\{t\geq0:\ \mbox{a}\ \phi_N\mbox{-event occurs at }t\ \}.
  \eeqnn
\blemma\label{lemma 2}
$\sigma^N_0\overset{p}{\longrightarrow}0$ as $N\rightarrow\infty$.
\elemma

\proof
By Lemma 1.3, we have for any $T$ and $0<\delta<(\theta/\gamma)^{1/\alpha}$, as $N\rightarrow\infty$,
  \beqlb
 \bfP\Big(\sup_{0\leq t\leq T}|\tilde{\eta}_N(t)-(\theta/\gamma)^{1/\alpha}|>\delta\Big)\rightarrow0.
 \label{constant}
  \eeqlb
Fix above $\delta$. Let $c_1=(\theta/\gamma)^{1/\alpha}-\delta$ and $c_2=(\theta/\gamma)^{1/\alpha}+\delta$.
Conditioned on $(\bfY_N(t),\tilde{\eta}_N(t))=(\bfy,u)$
with $(\bfy,u)\in(E\setminus\Pi)\times[c_1,c_2]$ at the current time $t$,
 \beqnn
 \bfP(\mbox{the
next event is a } \phi^N\mbox{-event} )=\frac{c_{\phi}^N(\bfy,u)}
 {c_{\psi}^N(\bfy)+c_{\phi}^N(\bfy,u)+c_{\Gamma}^N(\bfy,u)}\leq C\varepsilon_{_N},
 \eeqnn
for some positive constant $C$;
\beqnn
 \bfP(\mbox{the
next event is a } \psi^N\mbox{-event} )=\frac{c_{\psi}^N(\bfy)}
 {c_{\psi}^N(\bfy)+c_{\phi}^N(\bfy,u)+c_{\Gamma}^N(\bfy,u)}\leq
 \frac{2n^3}{2n^3+\theta+\gamma c_1^\alpha},
 \eeqnn
for sufficiently large $N$;
 \beqnn
 \bfP(\mbox{the
next event is a } \Gamma^N\mbox{-event} )=\frac{c_{\Gamma}^N(\bfy,u)}
 {c_{\psi}^N(\bfy)+c_{\phi}^N(\bfy,u)+c_{\Gamma}^N(\bfy,y)}\leq
 \frac{\theta+\gamma c_2^\alpha}{2+\theta+\gamma c_2^\alpha},
 \eeqnn
for sufficiently large $N$. Inspired by Taylor and V\'{e}ber \cite[Lemma 3.1]{TV09}, we fix
some $s>0$ and consider
 $$
 \bfP(\sigma_0^N>s)= \bfP(D)+o(1),
 $$
where
$$
 D = \{\sigma_0^N>s,\ \sup_{0\leq t\leq s}|\tilde{\eta}_N(t)-(\theta/\gamma)^{1/\alpha}|\leq\delta\}.
 $$
Then
 \beqnn
 \bfP(D)
 \ar=\ar \bfP(\{\mbox{at most } n\ \psi^N\mbox{-events occur in} [0,s]\}\cap D)\\
 \ar=\ar\bfP(\{\mbox{at most } n\ \psi^N\mbox{-}\mbox{ and  at least a } \phi^N\mbox{-events occur in} [0,s]\}\cap D)\\
 \ar \ar
 +\,\bfP(\{\mbox{at most } n\ \psi^N\mbox{-}\mbox{ and  no } \phi^N\mbox{-events occur in} [0,s]\}\cap D)\\
 \ar=:\ar I_1+I_2.
 \eeqnn
Note that we have $Y^N(t)\in E\setminus\Pi$ for $t\in[0,s]$ if $\sigma^N_0>s$.
Let $p=\frac{2n^3}{2n^3+\theta+\gamma c_1^\alpha}\vee
\frac{\theta+\gamma c_2^\alpha}{2+\theta+\gamma c_2^\alpha}$.
Then
 \beqlb
 I_1\ar\leq\ar\sum_{k=0}^n \bfP(
 \{\mbox{exactly } k\ \psi^N\mbox{-events}
 \mbox{ before a }
 \phi^N\mbox{-event occur in } [0,s]\}\cap D)\nonumber\\
 \ar=\ar
 \sum_{k=0}^n\sum_{l=0}^\infty\bfP(
 \{\mbox{exactly } k\ \psi^N\mbox{- and }l\ \Gamma^N\mbox{-events }\mbox{before a }
 \phi^N\mbox{-event in } [0,s]\}\cap D)
 \nonumber\\
 \ar\leq\ar
 \sum_{k=0}^n\sum_{l=0}^\infty\binom{k+l}{k}p^{k+l}
 (C\varepsilon_{_N}),
 \label{I_1}
 \eeqlb
Since $0<p<1$, $\sum_{k=0}^n\sum_{l=0}^\infty\binom{k+l}{k}p^{k+l}<\infty$. Then $I_1\rightarrow0$ as $N\rightarrow\infty$.
Let $U_j^N$ be the arrival time
 of the $j$'th event occurring to $(\bfY_N,\tilde{\eta}_N)$. For $I_2$,

 \beqnn
 I_2
 \ar=\ar
\sum_{k=0}^n\sum_{l=0}^\infty\bfP(
 \{\mbox{exactly } k\ \psi^N\mbox{-events, }l\ \Gamma^N\mbox{-events }\mbox{and no } \phi^N\mbox{-events occur in } [0,s]\}\cap D)\\
 \ar\leq\ar
\sum_{k=0}^n\sum_{l=0}^\infty\binom{k+l}{k}p^{k+l}\bfP(\{U_{k+l}^N<s, U^N_{k+l+1}>s\}\cap D).
 \eeqnn
Conditioned on $(\bfY_N(t),\tilde{\eta}_N(t))=(\bfy,u)$
with $\bfy\in(E\setminus\Pi)$,
the rate for the event occurring to $(\bfY_N,\tilde{\eta}_N)$
at time $t$ is $c_{\psi}^N(\bfy)+c_{\phi}^N(\bfy,u)+c_{\Gamma}^N(\bfy,u)$ and
$c^N_{\psi}(\bfy)\geq 2/\varepsilon_{_N}$. Then $U^N_{k+l+1}$ is stochastically bounded by
the sum of $k+l+1$ i.i.d. exponential variables with parameter $2/\varepsilon_{_N}$
whose distribution becomes concentrated close to $0$ as $N\rightarrow\infty$.
Thus as $N\rightarrow\infty$, $\bfP(\{U_{k+l}^N<s, U^N_{k+l+1}>s\}\cap D)\rightarrow0$
and by the dominated convergence theorem, $I_2\rightarrow0$.\qed

\blemma\label{lemma 3} There exist positive constants $M$ and $K_1$ such that for any $s>0$,
 \beqnn
 \limsup_{N\rightarrow\infty}\bfP(\sigma^N_1\leq s)\leq M(1-e^{-K_1 s}).
 \eeqnn
\elemma

\proof By the proof of $(\ref{I_1})$,
$\bfP(\mbox{at least one }\phi_{N}\mbox{-event occurs before } \sigma_0^N)\rightarrow0$
as $N\rightarrow\infty$.
Then by (\ref{constant}), we have
$$
\bfP(\sigma_1^N\leq s) = \bfP(G)+o(1),
$$
where
 $$
 G=
\{\mbox{only }\psi_N\mbox{- or }\Gamma_N\mbox{- events before } \sigma_0^N,
 \sup_{0\leq t\leq s}|\tilde{\eta}_N(t)-(\theta/\gamma)^{1/\alpha}|
 \leq\delta\mbox{ and } \sigma_1^N\leq s\}.
$$
Recall that $\bfY_N(0)=\bfy$. If only $\psi_N$- or $\Gamma_N$- events occur before $\sigma_0^N$,
$\sigma_0^N<\sigma_1^N$ and $\bfY_N(t)=\bar{\bfy}$ for $t\in[\sigma_0^N,\sigma_1^N]$.
Furthermore $Y_N(t)\in\Pi$ and $c_{\psi}^N\equiv0$ for $t\geq\sigma_0^N$.
Conditioned on $(\bfY_N(t),\tilde{\eta}_N(t))=(\bfy,u)$
with $(\bfy,u)\in\Pi\times[c_1,c_2]$,
\beqnn
K_1\leq c_{\phi}^N(\bfy,u)\leq K_2,\quad
\varepsilon^{-1}_{_N}(\theta+\gamma c_1^{\alpha})/2\leq c_{\Gamma}^N(\bfy,u)\leq\varepsilon_{_N}^{-1}(\theta+\gamma c_2^\alpha),
\eeqnn
for sufficiently large $N$, where $K_1=[\gamma(c_1^{\alpha-1}\wedge c_2^{\alpha-1})]\wedge\theta$
and $K_2=n^2+n\theta+\gamma n(c_1^{\alpha-1}\vee c_2^{\alpha-1})$. Then
 \beqnn
 \frac{c_{\phi}^N(\bfy,u)}
 {c_{\phi}^N(\bfy,u)+c_{\Gamma}^N(\bfy,u)}
 \leq
\frac{2K_2\varepsilon_{_N}}
{2K_2\varepsilon_{_N}+\theta+\gamma c_1^\alpha},\quad
\frac{c_{\Gamma}^N(\bfy,u)}
 {c_{\phi}^N(\bfy,u)+c_{\Gamma}^N(\bfy,u)}
 \leq
\frac{\theta+\gamma c_2^\alpha}
{K_1\varepsilon_{_N}+\theta+\gamma c_2^\alpha}.
 \eeqnn
For $(\bfY_N(\cdot), \tilde{\eta}_N(\cdot))$ with initial value $(\bfy,u)\in\Pi\times[c_1,c_2]$, recall that
$U_j^N$ denotes the arrival time of the $j$'th event occurring to $(\bfY_N,\tilde{\eta}_N)$ and $U_0^N=0$.
It is not hard to see that $U^N_j$ is stochastically larger than
the sum of $j$ i.i.d. exponential variables with parameter
$\varepsilon_{_N}^{-1}(\theta+\gamma c_2^\alpha)+K_2$. We have
 \beqnn
  \bfP(G)
 \ar=\ar
\sum_{k=0}^\infty
\bfP(\{\mbox{exactly } k\ \Gamma_N\mbox{-events occur in } [\sigma_0^N,\sigma_1^N]\}\cap G)\\
 \ar\leq\ar
 \sum_{k=0}^\infty
\Big(\frac{\theta+\gamma c_2^\alpha}
{K_1\varepsilon_{_N}+\theta+\gamma c_2^\alpha}\Big)^k
\frac{2K_2\varepsilon_{_N}}
{2K_2\varepsilon_{_N}+\theta+\gamma c_1^\alpha}
\bfP\Big(\sigma_0^N+U_k^N\leq s\Big)\\
 \ar\leq\ar
 \sum_{k=0}^\infty
\Big(\frac{\theta+\gamma c_2^\alpha}
{K_1\varepsilon_{_N}+\theta+\gamma c_2^\alpha}\Big)^k
\frac{2K_2\varepsilon_{_N}}
{2K_2\varepsilon_{_N}+\theta+\gamma c_1^\alpha}
\bfP\Big(\sigma_0^N+\sum_{j=1}^k \tilde{V}^N_j\leq s\Big)\\
 \ar\leq\ar
M\bfP\Big(\sigma_0^N+\sum_{j=1}^{T_N}\tilde{V}_j^N\leq s\Big),
 \eeqnn
for some positive constant $M$ and sufficiently large $N$, where $\{\tilde{V}_j^N\}$ are i.i.d. exponential variables with parameter
$\varepsilon_{_N}^{-1}(\theta+\gamma c_2^\alpha)+K_2$, and $T_N$ is a geometric variable
with parameter $\frac{K_1\varepsilon_{_N}}{K_1\varepsilon_{_N}+\theta+\gamma c_2^\alpha}$ independent of
$\{\tilde{V}_j^N\}$. Since $\sigma_0^N\overset{p}{\rightarrow}0$, a simple calculation shows that
$\sigma_0^N+\sum_{j=0}^{M_N}V_j^N$ converges
weakly to an exponential variable with parameter $K_1$. The lemma is proved.\qed

\blemma\label{lemma 4} \;Under conditions (A) and (B),
the ancestral process
$\{\bfY_N(t),\ t\geq0\}$ starting at $\bfy$ with $\bfy\in\Pi$ converges weakly on $D([0,\infty), \Pi)$ to
$\{\bfY(t),\ t\geq0\}$ given by (\ref{Y}) starting at ${\bfy}$.
\elemma

\proof If the process $Y_N(t)$ stays
in the space of $\Pi$, $\psi_N$ and the fourth term in $\phi_N$ vanishes.
Then for any bounded function $g$ on $E\times\mathbb{R}_+$
define $B_Ng=\tilde{\phi}_Ng+\Gamma_Ng$,
where $\Gamma_N$ is given in (\ref{A_N}) and
  \beqnn
  \tilde{\phi}_Ng(\bfy,u)
  \ar=\ar
  2\binom{y_0}{2}\Big(g(\bfy-\bfe_0,u)-g(\bfy,u)\Big)\\
  \ar \ar
  +\,\theta y_0
  \Big(g(\bfy-\bfe_0+\bfe_1,u+\varepsilon_{_N})-g(\bfy,u)\Big)\\
  \ar \ar
  +\,\gamma u^{\alpha-1}y_1
  \Big(g(\bfy-\bfe_1+\bfe_0,u-\varepsilon_{_N})-g(\bfy,u)\Big)1_{\{u>0\}}\\
  \ar \ar
  +\,\Big(\varepsilon_{_N}R_{1,N}g(\bfy,u)+\frac{1}{N}u^{\alpha-1}1_{\{u>0\}}
  R_{2,N}g(\bfy,u)\Big).
    \eeqnn Let $\mathcal{F}^N_t=
\{(\bfY_N(s),\tilde{\eta}_N(s)):\ 0\leq s\leq t\}$. Because of the Markov property of $(\bfY_N(\cdot),\tilde{\eta}_N(\cdot))$,
 \beqnn
 g(\bfY_N(t),\tilde{\eta}_N(t))-g(\bfy,\tilde{\eta}_N(0))-\int_0^t(B_Ng)(\bfY_N(s),\tilde{\eta}_N(s))ds
 \eeqnn
is a local $(\mathcal{F}^N_t)$-martingale. Let
 \beqnn
 \tau_N=\inf\{
t\geq0,\ |\tilde{\eta}_N(t-)-(\theta/\gamma)^{1/\alpha}|\geq\delta\mbox{ or } |\tilde{\eta}_N(t)-(\theta/\gamma)^{1/\alpha}|\geq\delta
 \},
 \eeqnn
where $\delta$ is a positive constant satisfying $0<\delta<(\theta/\gamma)^{1/\alpha}$. For sufficiently large $N$,
there exist $c_1>0$ and $c_2>0$ such that
$c_1\leq \tilde{\eta}_N(t\wedge \tau_N)\leq c_2$ for any $t\geq0$.
Then for any function $f$ on $E$,
 \beqnn
 \zeta_N(t)=f(\bfY_N(t\wedge\tau_N))-f(\bfy)
 -\int_0^{t\wedge\tau_N}(B_Nf)(\bfY_N(s\wedge\tau_N),\tilde{\eta}_N(s\wedge\tau_N))ds
 \eeqnn
is a bounded martingale. Indeed, for sufficiently large $N$,
$|(B_Nf)(\bfY_N(s\wedge\tau_N),\tilde{\eta}_N(s\wedge\tau_N))|\leq
2\big(n^2+n\theta+n\gamma(c_1^{\alpha-1}\vee c_2^{\alpha-1})\big)\|f\|$
for any $s\geq0$.
Note that $E$ is a finite set, so the discrete topology on $E$
makes it a complete and compact metric space and
any real valued function $f$ on $E$ is bounded and continuous.
By Ethier and Kurtz \cite[Theorem 9.1 and 9.4, p.142]{EK86},
the process $\{Y_N(t\wedge\tau_N)\}$ is relatively compact. On the other
hand, for any $t\geq0$,
 \beqnn
 \bfP(\tau_N<t)\leq\bfP\Big(\sup_{0\leq s\leq t}|\tilde{\eta}_N(s)-(\theta/\gamma)^{1/\alpha}|\geq\delta\Big)\rightarrow0,
 \eeqnn
which shows that $\tau_N\overset{p}{\rightarrow}\infty$.
Let $\{\bfY(t)\}$ be any limit point of $\{\bfY_N(t\wedge\tau_N)\}$.
By Skorokhod's representation theorem we may assume that on some Skorokhod
space $(\tau_N,\bfY_N(t\wedge\tau_N),\tilde{\eta}_N(t\wedge\tau_N),)
\overset{a.s.}{\rightarrow}(\infty,\bfY(t),(\theta/\gamma)^{1/\alpha})$ in
the topology of $\mathbb{R}_+\times D([0,\infty),E\times\mathbb{R}_+)$. Thus
$\zeta_N(t)\overset{a.s.}{\rightarrow}\zeta(t)$ and $\zeta(t)$ is given by
 \beqlb
 \zeta(t)=f(\bfY(t))-f(\bfx)
 -\int_0^{t}(Bf)(\bfY(s))ds,\label{martingale problem}
 \eeqlb
where
 \beqnn
 Bf(u,\bfy)\ar=\ar2\binom{y_0}{2}\Big(f(\bfy-\bfe_0)-f(\bfy)\Big)
  +\theta y_0
  \Big(f(\bfy-\bfe_0+\bfe_{1})-f(\bfy)\Big)\\
  \ar \ar
  +\,\theta(\gamma/\theta)^{1/\alpha}y_1
  \Big(f(\bfy-\bfe_1+\bfe_0)-f(\bfy)\Big).
  \eeqnn
Since $\sup_N|\zeta_N(t)|<\infty$, $\zeta_N(t)\overset{L_1}{\rightarrow}\zeta(t)$
for any $t\geq0$. Thus $\zeta(t)$ is a martingale. Since $\{\bfY(t)\}$ given by (\ref{Y}) is
the unique solution to the martingale problem (\ref{martingale problem}),
We have that $\{\bfY_N(t\wedge\tau_N)\}$ converges weakly to $\{\bfY(t)\}$ given by (\ref{Y})
in $D([0,\infty),E)$. Furthermore, for any $\epsilon>0$ and any $t\geq0$,
 \beqnn
\bfP\Big(\sup_{0\leq s\leq t} |Y_N(s\wedge\tau_N)-Y_N(s)|>\epsilon\Big)
\leq\bfP(\tau^N<t)\rightarrow0,
 \eeqnn
as $N\rightarrow\infty$. The lemma follows from the above limit.\qed

{\it Proof of Theorem \ref{thm}\;} Let $\bfP_{\bfy}(\cdot)$ be the distribution
of $(\bfY_N(\cdot),\tilde{\eta}_N(\cdot))$ with initial value $(\bfy,\tilde{\eta}_N(0))$, where
$\tilde{\eta}_N(\cdot)$ is distributed as $\pi_N$ given in Section 2.
Let $f_1,\cdots, f_k$ be real-valued functions on $E$. Choose
$0<s<t_1<\cdots<t_k$. Let $Q_N=\{\sigma_0^N<s<\sigma_1^N\}$.
Then
 \beqnn
 \ar\ar\bfE_{\bfy}\Big[\Pi_{i=1}^kf_i(\bfY_N(t_i))1_{Q_N}\Big]\\
 \ar=\ar\bfE_{\bfy}\Big[1_{Q_N}\bfE\Big[\Pi_{i=1}^kf_i(\bfY_N(t_i))|\mathcal{F}^N_s\Big]\Big]\\
 \ar=\ar
\bfE_{\bfy}\Big[1_{Q_N}\bfE_{(\bar{\bfy},\,\tilde{\eta}_N(s))}\Big[\Pi_{i=1}^kf_i(\bfY_N(t_i-s))\Big]\Big]\\
 \ar=\ar
\bfE_{\bfy}\Big[\bfE_{(\bar{\bfy},\,\tilde{\eta}_N(s))}\Big[\Pi_{i=1}^kf_i(\bfY_N(t_i-s))\Big]\Big]
-\bfE_{\bfy}\Big[1_{\bar{Q}_N}\bfE_{((\bar{\bfy},\,\tilde{\eta}_N(s))}\Big[\Pi_{i=1}^kf_i(\bfY_N(t_i-s))\Big]\Big]\\
\ar=\ar
\bfE_{\bar{\bfy}}\Big[\Pi_{i=1}^kf_i(\bfY_N(t_i-s))\Big]
-\bfE_{\bfy}\Big[1_{\bar{Q}_N}\bfE_{(\bar{\bfy},\,\tilde{\eta}_N(s))}\Big[\Pi_{i=1}^kf_i(\bfY_N(t_i-s))\Big]\Big].
 \eeqnn
where $\bar{Q}_N$ is the complement of the set $Q_N$. The last equality follows from the fact that $\tilde{\eta}_N(\cdot)$ is stationary. Then by Lemmas \ref{lemma 2},
\ref{lemma 3} and \ref{lemma 4},
 \beqnn
\ar \ar\limsup_{N\rightarrow\infty}
\bigg|\bfE_{\bfy}\Big[\Pi_{i=1}^kf_i(\bfY_N(t_i))\Big]
-\bfE_{\bar{\bfy}}\Big[\Pi_{i=1}^kf_i(\bfY(t_i-s))\Big]\bigg|\\
\ar\leq\ar
 2\max_i\|f_i\|\limsup_{N\rightarrow\infty}\bfP(\bar{Q}_N)+
\limsup_{N\rightarrow\infty}
\bigg|\bfE_{\bar{\bfy}}\Big[\Pi_{i=1}^kf_i(\bfY_N(t_i-s))\Big]
-\bfE_{\bar{\bfy}}\Big[\Pi_{i=1}^kf_i(\bfY(t_i-s))\Big]\bigg|\\
\ar\leq\ar
2M\max_i\|f_i\|(1-e^{-K_1 s}),
 \eeqnn
which goes to $0$ as $s\rightarrow0$. Note that $\bfY(t)$ is stochastically
continuous. $\bfE_{\bar{\bfy}}\Big[\Pi_{i=1}^kf_i(\bfY(t_i-s))\Big]$ converges to
$\bfE_{\bar{\bfy}}\Big[\Pi_{i=1}^kf_i(\bfY(t_i))\Big]$ as $s\rightarrow0$.
Then we have that $\lim_{N\rightarrow\infty}
\bfE_{\bfy}\Big[\Pi_{i=1}^kf_i(\bfY^N(t_i))\Big]=
\bfE_{\bar{\bfy}}\Big[\Pi_{i=1}^kf_i(\bfY(t_i))\Big]$.\qed

{\it Proof of Theorem \ref{thm2}\;} Step 1: recall the notation in Section 2.2. Under the homeomorphism,
the subspace $\Pi$ can be regarded as $\Gamma(\Pi)$ for simplicity.
 It follows from (\ref{martingale problem}) that for any function $f$ on $\Pi$,
 \beqnn
f(\bfY_k(t))-f(\bfy)
 -\int_0^{t}(B_kf)(\bfY_k(s))ds
 \
 \eeqnn
is a martingale, where
 \beqnn
 B_kf(\bfy)
  \ar=\ar
 2\binom{y_0}{2}\Big(f(\bfy+(-1,0))-f(\bfy)\Big)
  +\theta_k y_0
  \Big(f(\bfy+(-1,1))-f(\bfy)\Big)\\
  \ar \ar
  +\,\theta_k(\gamma_k/\theta_k)^{1/\alpha}y_1
  \Big(f(\bfy+(1,-1))-f(\bfy)\Big).
  \eeqnn
Recall that $\bfY_k(t)=(Y_k^0(t),Y_k^1(t))$ and $Y_k(t)=Y_k^0(t)+Y_k^1(t)$. For any function $g$
on $I_n$, let $f(\bfy)=g(y_0+y_1)$ for $\bfy\in\Pi$. Then
\beqnn
g(Y_k(t))-g(n)
 -\int_0^{t}(\tilde{B}_kg)(Y_k^0(s),Y_k^1(s))ds
 \eeqnn
is also a martingale, where $\tilde{B}_kg(y)=2\binom{y_0}{2}(g(y-1)-g(y))$. Note that $I_n$
is a finite set, so the discrete topology on $I_n$ makes it a complete and compact metric space.
Any real valued function $g$ on $I_n$ is bounded and continuous. Then $Y_k(\cdot)$ satisfies
the compact containment condition. For each $T>0$,
$\sup_k\int_0^T|\tilde{B}_kg(Y_k^0(s),Y_k(s))|ds\leq 2n^2T\|g\|$, where $\|g\|=\sup_{y\in I_n}
|g(y)|$. By Ethier and Kurtz \cite[Theorem 9.1 and 9.4, p.142]{EK86},
$Y_k(\cdot)$ is relatively compact in $D([0,\infty), I_n)$.

Step 2: suppose that $\{\xi^k_j(\cdot)\}_{j=1}^n$
is the sequence of i.i.d.\,Markov chains taking values in $\{0,1\}$ and whose transition rate matrix is given by
 \beqnn
 \left( \begin{array}{ll}
     -1 &
     1\\
     (\frac{\gamma_k}{\theta_k})^{\frac{1}{\alpha}} & -(\frac{\gamma_k}{\theta_k})^{\frac{1}{\alpha}}   \\
        \end{array} \right).
 \eeqnn
Let $P^k_{ij}(t)=\bfP(\xi^k_1(t)=j|\xi^k_1(0)=i)$. A simple calculation shows that
 \beqnn
 P_{00}^k(t)=1-P^k_{01}(t)\ar=\ar
\frac{(\theta_k/\gamma_k)^{1/\alpha}}{1+(\theta_k/\gamma_k)^{1/\alpha}}
+\frac{1}{1+(\theta_k/\gamma_k)^{1/\alpha}}e^{-(1+(\theta_k/\gamma_k)^{1/\alpha})t},\\
 P_{10}^k(t)=1-P_{11}^k(t)\ar=\ar
\frac{(\theta_k/\gamma_k)^{1/\alpha}}{1+(\theta_k/\gamma_k)^{1/\alpha}}
-\frac{(\theta_k/\gamma_k)^{1/\alpha}}{1+(\theta_k/\gamma_k)^{1/\alpha}}e^{-(1+(\theta_k/\gamma_k)^{1/\alpha})t}.
 \eeqnn
Let $\zeta_n^k(t)=\sum_{j=1}^n1_{\{\xi_j^k(t)=0\}}$.
Since $\{\xi_i^k(t)\}_{i=1}^n$ are independent of each other, it is not hard to see that
for any $g$ on $I_n$,
 \beqnn
 \sup_{x,y \in I_n}\Big|\bfE_x[g(\zeta_n^k(t))]
 -\bfE_y[g(\zeta_n^k(t))]\Big|\leq2n\|g\|e^{-(1+(\theta_k/\gamma_k)^{1/\alpha})t},\ t\geq0.
 \eeqnn
This implies $\zeta_n^k(t)$ satisfies the $\phi$-mixing condition (see \cite[P.111]{B05}). By (1.13) of \cite[p.109]{B05},
 \beqnn
 \ar \ar \sup_{y\in I_n}\big|\bfE_y[g(\zeta_n^k(t_2))g(\zeta_n^k(t_1))]-\bfE_y[g(\zeta_n^k(t_2))]\bfE_y[g(\zeta_n^k(t_1))]\big|\\
  \ar\leq\ar 2\sqrt{2n}\|g\|^2e^{-(1+(\theta_k/\gamma_k)^{1/\alpha})(t_2-t_1)/2}
 \eeqnn
for any $t_2\geq t_1\geq0$. Then
 \beqlb
 \ar\ar\bfE_y\Big[\Big(\int_0^t\Big(g(\zeta_n^k(\theta_ks))-\bfE_y[g(\zeta_n^k(\theta_ks))]\Big)ds\Big)^2\Big]\nonumber\\
 \ar=\ar\bfE_y\bigg[\int_0^t\int_0^tds_1ds_2\Big(g(\zeta_n^k(\theta_ks_1))-\bfE_y[g(\zeta_n^k(\theta_ks_1))]\Big)
  \Big(g(\zeta_n^k(\theta_ks_2))-\bfE_y[g(\zeta_n^k(\theta_ks_2))]\Big)\bigg]\nonumber\\
  \ar=\ar \int_0^t\int_0^tds_1ds_2
  \Big(\bfE_y[g(\zeta_n^k(\theta_ks_2))
  g(\zeta_n^k(\theta_ks_1))]-\bfE_y[g(\zeta_n^k(\theta_ks_2))]\bfE_y[g(\zeta_n^k(\theta_ks_1))]\Big)\nonumber\\
  \ar\leq\ar
  C(n)\|g\|^2\int_0^tds_2\int_0^{t}e^{-\theta_k(1+(\theta_k/\gamma_k)^{1/\alpha})|s_2-s_1|/2}ds_1\nonumber\\
  \ar\leq\ar C(n)\|g\|^2t/\theta_k,
   \label{inequility}
 \eeqlb
where $C(n)$ is a constant only depending $n$.  Since $P_{00}^k(\theta_kt)\rightarrow\frac{p^{1/\alpha}}{1+p^{1/\alpha}}$
and $P_{01}^k(\theta_kt)\rightarrow\frac{1}{1+p^{1/\alpha}}$ as $k\rightarrow\infty$, it is easy to see for any $t\geq0$,
$\zeta_n^k(\theta_kt)\overset{d}{\rightarrow}\zeta_n$ as $k\rightarrow\infty$, where
$\zeta_n$ follows the Binomial distribution, i.e., $\zeta_n\sim Bn(n,\frac{p^{1/\alpha}}{1+p^{1/\alpha}})$. Note that $I_n$ is finite.
The dominated convergence theorem shows that
 \beqnn
 \sup_{y\in I_n}\int_0^t\Big|\bfE_y[g(\zeta_n^k(\theta_ks))]-\bfE[g(\zeta_n)]\Big|ds\rightarrow0,
 \eeqnn
as $k\rightarrow\infty$. Combined with (\ref{inequility}), we have as $k\rightarrow\infty$,
 \beqlb
 \sup_{y\in I_n}\bfE_y\Big[\Big(\int_0^t\Big(g(\zeta_n^k(\theta_ks))-\bfE[g(\zeta_n)]\Big)ds\Big)^2\Big]\rightarrow0.
 \label{convergence}
  \eeqlb

Step 3: $(Y^0_k(t),Y_k(t))$ is a Markov process as in Step 1. $Y^0_k(0)=y\in I_n$ and $Y_k(0)=n$.
Let $\mathcal{F}^k_t=\sigma\{(Y^0_k(s),Y_k(s)):\ 0\leq s\leq t\}$. Define
$T^k_j=\inf\{t\geq0: Y_k(t)=n-j\}$ with $T^k_0=0$
and $\tau_j^k=T_j^k-T_{j-1}^k$. Set $h(y)=y(y-1)$
for $y\in I_n$.
By (\ref{convergence}), we have
 \beqnn
\bfP(\tau^k_{j+1}>t)\ar=\ar\bfE\Big[\bfP\Big(\tau^k_{j+1}>t\Big|\mathcal{F}^k_{T^k_{j}}\Big)\Big]\nonumber\\
 \ar=\ar\bfE\Big[\bfP_{(Y_k^0(T^k_{j}),n-j)}(\tau^k_{j+1}>t)\Big]\nonumber\\
 \ar=\ar\bfE\Big[\bfE_{Y_k^0(T^k_{j})}\Big(\exp\Big\{-\int_0^th(\zeta^k_{n-j}(\theta_ks))ds\Big\}\Big)\Big]\nonumber\\
 \ar\rightarrow\ar e^{-\bfE[h(\zeta_{n-j})]t},
  \label{convergence2}
  \eeqnn
as $k\rightarrow\infty$. Similarly
 \beqnn
\bfP(\tau^k_1>t,\tau^k_2>s)=
\bfE\Big[1_{\{\tau^k_1>t\}}\bfE_{Y_k^0(T^k_1)}\Big(\exp\Big\{-\int_0^th(\zeta^k_{n-1}(\theta_ks))ds\Big\}\Big)\Big].
 \eeqnn
Then
 \beqnn
 \ar\ar|\bfP(\tau^k_1>s,\tau^k_2>t)-e^{-\bfE[h(\zeta_n)]s-\bfE[h(\zeta_{n-1})]t}|\\
 \ar\leq\ar|\bfP(\tau^k_1>s)-e^{-\bfE[h(\zeta_n)]s}|+
 \sup_{y\in I_n}\bfE_y\Big[\Big|\int_0^t\Big(h(\zeta_{n-1}^k(\theta_ks))-\bfE[h(\zeta_{n-1})]\Big)ds\Big|\Big].
 \eeqnn
By (\ref{convergence}) and (\ref{convergence2}) we have that the second term in the right-hand side of the above inequality goes to $0$ as $k\rightarrow\infty$.
By induction, $(\tau_1^k,\cdots,\tau_{n-1}^k)\overset{d}{\rightarrow}(\tau_1,\cdots,\tau_{n-1})$, where $\{\tau_j\}_{j=1}^{n-1}$
is independent of each other and $\tau_j$ follows the exponential distribution with parameter $c_{n-j+1}$. It follows that
$\{Y_k(t),\,t\geq0\}$ converges in the sense of finite-dimensional distributions to the $n$-Kingman coalescent process $\{K(t),\,t\geq0\}$.
Since $\{Y_k(t)\}$ is relatively compact, the theorem is proved.\qed

\paragraph{Acknowledgments.} AL and CM were financially supported by grant MANEGE `Mod\`eles Al\'eatoires en \'Ecologie, G\'en\'etique et \'Evolution' 09-BLAN-0215 of ANR (French national research agency). 
AL also thanks the {\em Center for Interdisciplinary Research in Biology} (Coll\`ege de France) for funding. CM also thanks the financial support from the
National Natural Science Foundation of China (NSFC) (No.11001137 and No.11271204) and the China Scholarship Council (CSC).

\end{document}